\documentclass[reqno,11pt]{amsart}
\usepackage{amscd,amssymb}

\theoremstyle{plain}
\newtheorem{Thm}[subsection]{Theorem}
\newtheorem{Cor}[subsection]{Corollary}
\newtheorem{Lem}[subsection]{Lemma}
\newtheorem{Prop}[subsection]{Proposition}
\newtheorem{Conj}[subsection]{Conjecture}

\theoremstyle{definition}
\newtheorem{Def}[subsection]{Definition}

\theoremstyle{remark}

\newtheorem{Rem}[subsection]{Remark}

\newtheorem{conj}{Conjecture}

\numberwithin{equation}{section}
\renewcommand{\rm}{\normalshape}

\newif\ifShowLabels
\ShowLabelstrue
\newdimen\theight
\def\TeXref#1{%
    \leavevmode\vadjust{\setbox0=\hbox{{\tt
        \quad\quad  {\small \rm #1}}}%
    \theight=\ht0
    \advance\theight by \lineskip
    \kern -\theight \vbox to
    \theight{\rightline{\rlap{\box0}}%
    \vss}%
    }}%

\ShowLabelsfalse

\renewcommand{\sec}[2]{\section{#2}\label{S:#1}%
    \ifShowLabels \TeXref{{S:#1}} \fi}
\newcommand{\ssec}[2]{\subsection{#2}\label{SS:#1}%
    \ifShowLabels \TeXref{{SS:#1}} \fi}

\newcommand{\refs}[1]{Section ~\ref{S:#1}}
\newcommand{\refss}[1]{Section ~\ref{SS:#1}}

\newcommand{\reft}[1]{Theorem ~\ref{T:#1}}
\newcommand{\refl}[1]{Lemma ~\ref{L:#1}}
\newcommand{\refp}[1]{Proposition ~\ref{P:#1}}

\newcommand{\refe}[1]{\eqref{E:#1}}

\newenvironment{thm}[1]%
    { \begin{Thm} \label{T:#1}  \ifShowLabels \TeXref{T:#1} \fi }%
    { \end{Thm} }

\renewcommand{\th}[1]{\begin{thm}{#1} \sl }
\renewcommand{\eth}{\end{thm} }

\newenvironment{lemma}[1]%
    { \begin{Lem} \label{L:#1}  \ifShowLabels \TeXref{L:#1} \fi }%
    { \end{Lem} }
\newcommand{\lem}[1]{\begin{lemma}{#1} \sl}
\newcommand{\elem}{\end{lemma}}

\newenvironment{propos}[1]%
    { \begin{Prop} \label{P:#1}  \ifShowLabels \TeXref{P:#1} \fi }%
    { \end{Prop} }
\newcommand{\prop}[1]{\begin{propos}{#1}\sl }
\newcommand{\eprop}{\end{propos}}

\newenvironment{corol}[1]%
    { \begin{Cor} \label{C:#1}  \ifShowLabels \TeXref{C:#1} \fi }%
    { \end{Cor} }
\newcommand{\cor}[1]{\begin{corol}{#1} \sl }
\newcommand{\ecor}{\end{corol}}

\newenvironment{defeni}[1]%
    { \begin{Def} \label{D:#1}  \ifShowLabels \TeXref{D:#1} \fi }%
    { \end{Def} }
\newcommand{\defe}[1]{\begin{defeni}{#1} \sl }
\newcommand{\edefe}{\end{defeni}}

\newenvironment{remark}[1]%
    { \begin{Rem} \label{R:#1}  \ifShowLabels \TeXref{R:#1} \fi }%
    { \end{Rem} }
\newcommand{\rem}[1]{\begin{remark}{#1}}
\newcommand{\erem}{\end{remark}}

\newenvironment{conjec}[1]%
    { \begin{Conj} \label{Co:#1}  \ifShowLabels \TeXref{Co:#1} \fi }%
    { \end{Conj} }
\renewcommand{\conj}[1]{\begin{conjec}{#1} \sl }
\newcommand{\econj}{\end{conjec}}

\newcommand{\eq}[1]%
    { \ifShowLabels \TeXref{E:#1} \fi
       \begin{equation} \label{E:#1} }
\newcommand{\eeq}{ \end{equation} }

\newcommand{\prf}{ \begin{proof} }
\newcommand{\epr}{ \end{proof} }

\newcommand\alp{\alpha}

\newcommand\lam{\lambda}		\newcommand\Lam{\Lambda}
\newcommand\sig{\sigma}

\newcommand\calG{{\mathcal{G}}}

\newcommand\calK{{\mathcal{K}}}
\newcommand\calL{{\mathcal{L}}}
\newcommand\calM{{\mathcal{M}}}

\newcommand\calP{{\mathcal{P}}}

\newcommand\calV{{\mathcal{V}}}
\newcommand\calW{{\mathcal{W}}}

\newcommand\RR{\mathbb{R}}

\newcommand\ZZ{\mathbb{Z}}

\newcommand\CC{\mathbb{C}}

	\newcommand\gra{{\mathfrak{a}}}

	\newcommand\grg{{\mathfrak{g}}}

	\newcommand\grn{{\mathfrak{n}}}

	\newcommand\grt{{\mathfrak{t}}}

\newcommand\sdp{\times \hskip -0.3em {\raise 0.3ex
\hbox{$\scriptscriptstyle |$}}} 

\newcommand\Hom{\operatorname {Hom}}

\newcommand\og{{\overline{g}}}

\newcommand\oz{{\overline{z}}}

\newcommand\hata{{\widehat{a}}}
\newcommand\hatA{{\widehat{A}}}

\newcommand\hatg{{\widehat{g}}}

\newcommand\hatT{{\widehat{T}}}

\newcommand\hatW{{\widehat{W}}}

\newcommand\x{\times}
\newcommand\ten{\otimes}

\newcommand{\ra}{\rangle}
\newcommand{\la}{\langle}

\newcommand\Bun{\operatorname{Bun}}

\begin{document}

\title{Cartan decomposition for complex loop groups}
\author{Alexander Braverman and David Kazhdan}
\address{Department of Mathematics, Brown University, 151 Thayer street, Providencde RI, 02912, USA}
\address{Einstein Institute of Mathematics, The Hebrew University of Jerusalem}
\maketitle
\begin{abstract}
  In this note we formulate and prove a version of Cartan decomposition for holomorphic loop groups, similar to Cartan decomposition for $p$-adic loop groups, discussed in \cite{BK}, \cite{gar:ihes}. The main technical tool that we use is the (well-known) interpretation of twisted conjugacy classes in the holomorphic loop group in terms of principal holomorphic bundles on an elliptic curve.
\end{abstract}

\sec{}{Introduction and statement of the results}

\ssec{finite-int}{Complex groups}Let $G$ be a connected complex reductive group.
Let us choose a pair of opposite Borel subgroups $B, B^-$ with the maximal torus
$T=B\cap B^-$. We denote by $\Lam$ the coweight lattice of $T$.  We also denote by $\sig:G\to G$ the corresponding Cartan anti-involution (which acts as identity on $T$).

The group $G$ is also endowed with a complex conjugation
$g\mapsto \og$, whose group of fixed points is the split real group $G_{\RR}$. We may assume that $^-$ commutes with $\sig$ and
preserves $B$ and $B^-$. For any $g\in G$ let us set $g^*=\sig(\og)$.

We set $K$ to be the set of fixed points of the involution $g\mapsto (g^*)^{-1}$. This is a maximal compact subgroup of $G$.

 Let $T_K=T\cap K$. Let $\gra$ be the space of invariants
of $^*$ in $\grt=Lie(T)$ and let $A$ be its exponential. Obviously we have $T=T_K\cdot A$.
We say that an element $x\in \gra$ is dominant
 if $\la x,\alp\ra\geq 0$ for any positive root $\alp$. We denote by $\alp^+$ the set of dominant elements in $\gra$ and
 by $A^+$ its exponential. It is clear that $A^+$ is the set of $W$-orbits on $A$.

The classical Cartan decomposition says that the natural map $A^+\to K\backslash G/K$ is bijective. In other words, for any $g\in G$
there exists unique $a\in A^+$ such that $g$ and $a$ lie in the same double coset with respect to $K$.
In this note we would like to generalize this result to loop groups.

By a "representation" of $G$ we shall mean an algebraic finite-dimensional representation. If $\pi:G\to GL(V)$ is such a representation,
we say that $(\pi,V)$ is unitary if there exists a positive definite Hermitian form $(\cdot,\cdot)$ on $V$ such that
$(\pi(g)(v),w)=(v,\pi(g^*)(w)$ for every $g\in G$ and $v,w\in V$. It is well-known that every representation of $G$ has
a unitary structure.

\ssec{}{The loop group and its compact form}
Let $LG_{hol}$ denote the group of all holomorphic maps $\CC^*\to G$ and let $\calG_{hol}$ be semidirect product
of $LG_{hol}$ with $\CC^*$ (which acts on $LG_{hol}$ by loop rotation). Similarly, we can define the corresponding polynomial loop
group $LG_{pol}$ and the group $\calG_{pol}$.

For any $g(z)\in LG_{hol}$ the map $\og(\oz)$ is again in $LG_{hol}$ and we shall denote it by $\overline{g(z)}$. Define also an anti-involution
$\tau:LG_{hol}\to LG_{hol}$ by setting
$$
\tau(g(z))=g^{\sig}(z^{-1}).
$$
Both $\tau$ and $^-$ extend to $\calG_{hol}$ in a natural way (in particular, $\tau$ acts on the "loop rotation" $\CC^*$ as identity).
For any $g\in \calG_{hol}$ define $g^*=\tau(\og)$ (note that if $g\in G\subset LG_{hol}$ we get the same definition as before).
We set $\calK_{hol}$ to be the set of all $g\in \calG_{hol}$ such that $g^{-1}=g^*$.
This is the analog of $K$ in the loop case. Note that under $\varpi$ it projects to the unit circle
$S^1\subset \CC^*$.

We have the natural homomorphism $\varpi:\calG_{hol}\to \CC^*$. Define the semi-group
\eq{}
\calG_{hol}^+=\calK_{hol}\cup\{ g\in\calG_{hol}|\ |\varpi(g)|>1\}
\end{equation}
(note that $\varpi$ sends $\calK_{hol}$ to the unit circle).

All of the above definitions make sense when applied to the group $\calG_{pol}$.

\ssec{}{Cartan decomposition}
Let $\hatT=T\x \CC^*, \hatA=A\times \RR$.
Let $\hatW=W\rtimes \Lam$. This group naturally acts on $\hatA$ where $W$ just acts on the first factor
and $\lam\in\Lam$ acts by sending $(a,x)$ to $(a+x\lam,x)$.

Let also
\eq{}
\hatA^+=\{(a,q)\in \hatA| a\in A^+\ \text{and\ } \alp(a)\leq q\ \text{for any positive root $\alp$}\}
\end{equation}
It is well-known that $\hatA^+$ is a fundamental domain for the action of the group $\hatW$ on $\hatA$.

The main result of this note is the following affine Cartan decomposition
\th{main}
The
natural map
$$
\hatA^+\to\calK_{hol}\backslash\calG_{hol}^+/\calK_{hol}
$$
is a bijection.
\eth
\noindent
\ssec{}{Remarks}

1) We were informed by P.~Etingof that the same result was proved by him in 1990's (unpublished); apparently, his proof was quite different from ours.

2) It makes sense to ask whether the same statement holds for $\calG_{pol}$, i.e. whether
$\calG_{pol}^+$ is equal to $\calK_{pol}\hatA^+\calK_{pol}$. We strongly believe that it is NOT true, but we do not know a counterexample.

3) \reft{main} should be compared with the corresponding statements in \cite{gar:ihes} and \cite{BK}, where a similar statement
is proved when $LG_{hol}$ is replaced by the group $G(\textsf{k}[z,z^{-1}])$ where $\textsf{k}$ is a local non-archimedian field.
In addition in \cite{BK} this decomposition is used in order to define the {\em spherical Hecke algebra} for the central extension
of the group $G(\textsf{k}[z,z^{-1}])\rtimes \textsf{k}^*$. We do not know how to do this in our present context (which essentially corresponds to the case
$\textsf{k}=\CC$).

\ssec{}{Real variant}
There is a version of \reft{main} for real (as opposed to complex) loop groups. We shall leave the details to the reader (let us only
note that one needs to work with the group of real analytic loops).
\ssec{}{Organization of the paper}In \refs{fd} we gather some facts about the group $G$ that will be needed in the proof of \reft{main}.
In \refs{ell} we review the relation between conjugacy classes in $\calG_{hol}^+$ and semi-stable bundles on elliptic curves. \refs{proof}
is devoted to the actual proof.
\ssec{}{Acknowledgements}We would like to thank Pavel Etingof for helpful discussions.
Both authors were partially supported by the BSF.
The first author was partially supported by the NSF grant DMS-1200807 and by Simons Foundation.
The second author was also partially supported by the European Research Council grant 247049.

\sec{fd}{Some finite-dimensional statements}
In this Section we collect some statement about the finite-dimensional group $G$ that we are going to need for the proof of \reft{main}.
For the purposes of this section we assume that $G$ is just a reductive algebraic group over $\CC$ endowed with $^-$ and $\sig$  as in
\refss{finite-int} (i.e. we drop the assumption that $G$ is connected).
\ssec{}{The map $\eta$} Consider the map $\eta:G\to G$ defined by $\eta(g)=gg^*$. Clearly $\eta(g)^*=\eta(g)$ for any $g\in G$.
Our first statement is this:
\lem{imageofeta}
Let $x\in G$ be such that

a) $x=x^*$

b) There exists a faithful unitary representation $\pi:G\to GL(V)$ such that all eigen-values of $\pi(x)$ are positive
real numbers.

Then $x$ is $K$-conjugate to an element of $A^+$.
\elem
\prf
Clearly, it is sufficient to prove that is $K$-conjugate to an element of $A$. Let us choose a representation $V$ as above.
Then $\pi(x)$ is a self-adjoint positive-definite matrix. Hence $\pi(x)=\exp(X)$ for some $X\in gl(V)$, such that $X^*=X$ . Since $\pi(G)$ is a closed
algebraic subgroup of $GL(V)$, it follows that $X=d\pi(Y)$ with $Y^*=Y$ for some $Y\in \grg$ and hence $x=\exp(Y)$.
Now $Y$ is $K$-conjugate to an element of $\gra$ (cf. Lemma 2.1.9 in \cite{Wallach}), hence $x$ is $K$-conjugate to an element of $A$.
\epr
\cor{eta-fdim}
Assume that $x\in G$ satisfies properties a) and b) as in \refl{imageofeta}.
Then there exists $g\in G$ such that $x=\eta(g)$.
\ecor
\prf
By \refl{imageofeta} $x=kak^{-1}$ for some $k\in K, a\in A$. It is clear that $\eta|_A$ is a surjective map from $A$ to $A$.
Hence we can write $a=\eta(b)=bb^*$ for some $b\in A$. Hence
$$
x=kbb^*k^{-1}=(kbk)(kbk)^*=\eta(kbk).
$$
\epr

The second fact, that we are going to need in the proof of \reft{main} is this:
\lem{z}
Let $Z$ be a reductive subgroup of $G$ containing $T$ which is stable under $^-$ and $\sig$. Let $x\in Z$ be such that $x=\eta(h)$ for some $h\in G$.
Then $x=\eta(g)$ for some $g\in Z$.
\elem
\prf
We are going to apply \refl{imageofeta} to $Z$ instead of $G$. Clearly, $x$ satisfies condition a) of \refl{imageofeta}.
Condition b) is also clear for any representation $V$ of $Z$ which extends to a representation of $G$.
\epr
\sec{ell}{The group $\calG_{hol}$ and bundles on elliptic curves}
\ssec{}{From $\calG_{hol}$ to $G$-bundles}
Let $q\in \CC^*$ such that $|q|>1$.
Let $E_q$ denote the elliptic curve $\CC^*/q^{\ZZ}$. It is well-known that
isomorphism classes of $G$-bundles on $E_q$ are in one-to-one correspondence with conjugacy classes in $\pi^{-1}(q)$.
Moreover, under this correspondence the automorphism group of a $G$-bundle is naturally isomorphic to the centralizer
of the corresponding element inside $LG_{hol}$. Let us recall this correspondence.

First of all, a $G$-bundle on $E_q$ is the same as a $G$-bundle on $\CC^*$ which is endowed with a $\ZZ$-equivariant structure,
where the element $1\in \ZZ$ acts on $\CC^*$ by means of multiplication by $q$.
We can naturally identify $\pi^{-1}(q)$ with $LG_{hol}$. Then to any $g(z)\in LG_{hol}$ we associate the trivial
$G$-bundle on $\CC^*$, where the $\ZZ$-equivariant structure is equal to the standard one twisted by $g(z)$. We shall
denote the corresponding $G$-bundle on $E_q$ by $\calP_{g(z),q}=\calP_{\hatg}$ if we set $\hatg=(g(z),q)$.
It is easy to see that conjugating $\hatg$ by some $a(z)$
leads to an isomorphic $G$-bundle on $E_q$. Moreover, the centralizer of $\hatg$  inside $LG_{hol}$ is equal to the automorphism group of $\calP_{\hatg}$.

\ssec{}{Bundles of degree zero}
Fix $q$ as above and
let us denote by $\Bun_G$ the stack of $G$-bundles on $E_q$.
The connected components of $\Bun_G$ are in one-to-one correspondence with the elements of $\pi_1(G)$. We denote
by $c_{\Bun_G}:\Bun_G\to \pi_1(G)$ the corresponding map. We shall say that a $G$-bundle
$\calP$ has degree zero if it lies in the connected component of the trivial bundle (i.e. if $c_{\Bun_G}(\calP)=e$).
On the other hand, clearly we have a natural homomorphism $c_{LG}:LG_{hol}\to \pi_1(G)$ taking every $g(z)$ to the homotopy class
of $g|_{S^1}$. The following lemma is left to the reader:
\lem{con-comp}
For every $g(z)\in LG_{hol}$ we have $c_{Bun_G}(\calP_{g(z),q})=c_{LG}(g(z))$.
\elem
\ssec{}{Semi-stable bundles}
Recall that a vector bundle $\calM$ of degree zero
on an elliptic curve $E$ is called semi-stable if it has no sub-bundles of degree $>0$.
Also a principal $G$-bundle $\calP$ is called semi-stable of degree 0 if

1) $\calP$ has degree zero

2) For some (equivalently, for any) almost faithful representation $V$ of $G$
the induced vector bundle $\calP_V$ is semi-stable (automatically of degree zero because of 1))(cf. \cite{Balaji}, Def. 3.4, \cite{FM}, Theorem 2.2).

Also, we say that a $G$-bundle $\calP$ is split if the structure group of $\calP$ can be reduced to $T$. It is clear that if
$\calP$ is semi-stable of degree $0$ and split, then the corresponding $T$-bundle $\calP_T$ will have degree $0$ as well.

The folllowing theorem collects the knowledge about semi-stable bundles that we are going to need in this paper:
\th{semi-stable}
\begin{enumerate}
  \item
  Any bundle of the form $\calP_{g,q}$ where $g\in G$ is semi-stable of degree $0$. If $G=GL(n)$ then the converse is also true.
  \item
  A $G$-bundle $\calP$ is isomorphic to some $\calP_{t,q}$ with $t\in T$ if and only if $\calP$ is semi-stable of degree $0$ and split.
  \item
  Let $\calM$ be a semi-stable vector bundle on $E_q$ of degree zero. Then $\calM$ is split if and only if any short
  exact sequence
  $$
  0\to\calM_1\to\calM\to\calM_2\to 0
  $$
  where $\calM_1,\calM_2$ are vector bundles of degree $0$, splits.
  \item
  Let $\calP$ be a semi-stable $G$-bundle of degree $0$. Assume in addition that for some almost faithful representation $V$ of $G$ the vector bundle $\calP_V$ is split. Then $\calP$ is split (and thus by (2) it is isomorphic to a bundle of the form $\calP_{t,q}$ with $t\in T$).
  
\end{enumerate}
\eth
\prf
By definition statement (1) is true if and only if it is true for $G=GL(n)$. In this case it is proved by Atiyah \cite{Atiyah}.

For (2) let us note that by (1) any bundle of the form $\calP_{t,q}$ is semi-stable of degree $0$ by (1). On the other hand, since $t\in T$ its structure group can be reduced to $T$. For the converse statement it is enough to note that the set of (isomorphism classes of)
$T$-bundles of degree $0$ on $E_q$ is naturally isomorphic to $T/q^{\Lam}$, which precisely means that any $T$-bundle of degree zero is of the form $\calP_{t,q}$ for some $t\in T$.

(3) is proved in \cite{Atiyah}, so let us prove (4).

According to \cite{FMW}, Theorem 2.6 any semi-stable $G$-bundle $\calP$ of degree $0$ is $S$-equivalent to a $T$-bundle of degree zero
\footnote{In {\em loc. cit.} it is only proved for simply connected $G$, but under our "degree $0$" assumption the same proof goes through in general}. According to \cite{Schmitt}, Section 2.4 this means that any $\calP$ as above has a $B$-structure such that the induced
$T$-bundle $\calP_T$ is of degree $0$. Thus $\calP$ defines a class in $H^1(E_q, \calP_{T,\grn})$ where $\calP_{T,\grn}$ denotes the vector bundle associated with $\grn$ (as a representation of $T$). We would like to show that this class is $0$. Let 
us choose $(\pi,V)$ as in (3) and let $B_V$ be a Borel subgroup of $GL(V)$ which contains the image of $B$ under $\pi$ and let
$\grn_V$ denote the nilpotent radical of its Lie algebra. Then we get a natural map $\grn\to \grn_V$ of $T$-representations; moreover
$\grn$ is a direct summand of $\grn_V$. Hence it is enough to prove that the corresponding class in $H^1(E_q,\calP_{T,\grn_V})$ is equal to $0$. This means that the corresponding full flag of subbundles of $\calP_V$ splits, which follows from (3).

\epr
In the proof of \reft{main} we shall need the following result.
\prop{centralizer}
Let $q\in \CC^*$ such that $|q|\neq 1$ and let $a\in T$. Let $Z$ denote the centralizer of the element $(a,q)$ in $LG_{hol}$.
Then
\begin{enumerate}
  \item
  The natural map $ev_1:Z\to G$ (sending every $g(z)$ to $g(1)$ is a closed embedding whose image is a reductive subgroup of $G$ containing
  $T$.
  \item
  The image of $Z$ in $G$ is closed under $\sig$.
  \item If $(a,q)\in \hatA$, then $Z$ is stable under  $^-$.
\end{enumerate}
\eprop
\prf
The assertions (2) and (3) of \refp{centralizer} are clear, so let us prove the first assertion.
Is is also clear that $T\subset Z$ and that it is enough to prove (1) when $G=GL(n)$.
In this case we just need to show the following. Let $\calM$ is a split semi-stable vector bundle on $E_q$ of degree
$0$ and let $V$ be its fiber at $1\in E_q$. Then we must show that the natural map from $Aut(\calM)$ to $GL(V)$ is a closed
embedding and the image is reductive.

Now since the bundle $\calM$ is split and semi-stable of degree 0, it follows that there exist non-isomorphic
line bundles $\calL_1,...,\calL_r$ on $E_q$ and some finite-dimensional vector spaces $V_1,\cdots, V_r$ over $\CC$
such that
$$
\calM=\bigoplus\limits_{i=1}^r \calL_i\ten V_i.
$$
By choosing a trivialization of the fiber of each $\calL_i$ at $1\in E_q$ we can identify $V$ with
$\oplus_{i=1}^r V_i$. Since $\Hom(\calL_i,\calL_j)=0$ for $i\neq j$ it follows that the automorphism group of $\calM$
is  isomorphic to $\prod\limits_{i=1}^r GL(V_i)$ which is naturally a closed reductive subgroup of $GL(V)$.

\epr


\sec{proof}{Proof of \reft{main}}
\ssec{}{}
The proof of \reft{main} is based on the following
two results:

\lem{k-to-g}
Let $\hatg\in\calG_{hol}^+$, $\hata\in \hatA^+$. Then $\hatg$ and $\hata$ lie in the same $\calK_{hol}$-double coset if and only if
$\eta(\hatg)$ is conjugate to $\eta(\hata)$.
\elem

\lem{const}
Let $\hatg\in\calG_{hol}^+$. Then $\eta(\hatg)$ is conjugate to an element in $\hatA^+$.
\elem
\ssec{}{Proof of \refl{k-to-g}}
The "only if part" is clear. Indeed, for any $\hatg_1,\hatg_2\in \calG_{hol}$, if we have $\hatg_2=k \hatg_1 k'$ with $k,k'\in \calK_{hol}$
then
$$
\eta(\hatg_2)=k \hatg_1 k' (k')^*\hatg_1^* k^*=k\eta(\hatg_1)k^{-1}.
$$
Let us also note that if $\eta(\hatg_2)=k\eta(\hatg_1)k^{-1}$ with $k\in\calK_{hol}$, then the element $k'=k^{-1}\hatg_2\hatg_1^{-1}$ lies
in $\calK_{hol}$ and $\hatg_2=k \hatg_1 k'$.
Thus for the "only if" part it is enough to show that
if  $\eta(\hatg)=h\eta(\hata)h^{-1}$ for some $h\in \calG_{hol}$, then $h$ can be chosen to lie in $\calK_{hol}$. Clearly we may assume that
$h\in LG_{hol}$.
Note that since $\eta(\hatg)^*=\eta(\hatg)$ and $\eta(a)^*=\eta(a)$, we have
$$
h\eta(a)h^{-1}=(h^*)^{-1} \eta(a) h^*.
$$
Let $Z$ denote the centralizer of $\eta(a)$ in $LG_{hol}$. Then the above formula shows that
$h^*h\in Z$. We need to look for $x\in Z$ such that $hx=((hx)^*)^{-1}=(h^*)^{-1}(x^*)^{-1}$. In other words, we should
be looking for $x\in Z$ such that $x^*x=(hh^*)^{-1}$.
To summarize, we need to show that the element $(hh^*)^{-1}\in Z$ lies in the image of the map $\eta$ restricted to $Z$, which
is equivalent to showing that $hh^*$ lies in $\eta(Z)$.

Now, according to \refp{centralizer} we may think about $Z$ as of a reductive subgroup of $G$, containing $T$ and stable
under $\sig$ and $^-$ (this is done by identifying $Z$ with $ev_1(Z)$). Hence  to show that $\eta(h)\in \eta(Z)$ it is enough to check that $\eta(h)$ satisfies
the conditions a), b) of \refl{z}. Condition a) is completely clear and condition b) follows from the fact
$$
ev_1(\eta(h))=h(1)h^*(1)=h(1)h(1)^*\in \eta(G).
$$
\ssec{}{Proof of \refl{const}}

Let $\hatg=(g(z),q)$. Clearly we can assume that $q$ is real and positive. Then $\eta(\hatg)=(h(z),q^2)$, where $h(z)=g(z)g^*(q^{-1}\oz^{-1})$.
Note that we have
\eq{dual}
h(z)=h^*(q^{-1}\oz^{-1})
\end{equation}
Moreover, if $(\pi,V)$ is a a unitary representation of $G$ then for any $z\in \CC^*$ such that $|z|^2=q^{-1}$ (i.e. $z=q^{-1}\oz^{-1}$)
the matrix $\pi(h(z))$ is self-adjoint and positive-definite.

Let us first assume that $G=GL(n)$ and let $\calM$ be the vector bundle of rank $n$ associated with the principal bundle $\calP_{\eta(z)}$ on $E_{q^2}$. Then first of all $\calM$ has degree 0 (this immediately follows from \refl{con-comp}).

On the other hand, let $\iota:E_{q^2}\to E_{q^2}$ be induced by the map $z\mapsto q^{-1}\oz^{-1}$ on $\CC^*$ and let $\calM^*$ be the Hermitian dual of $\calM$. Then \refe{dual} implies that $\calM$ is isomorphic $\iota^*\calM^*$.
Moreover, let
$$
S_q=\{ z\in\CC^*|\ |z|^2=q^{-1}\}.
$$
We can (and will) identify $S_q$ with its image in $E_{q^2}$. Then the above map $\calM\to \iota^*\calM^*$ defines a Hermitian pairing
on $\calM|_{S_q}$ and it follows from the remark after \refe{dual} that this pairing is positive-definite. Hence for any subbundle
$\calV\subset \calM$ the composition map
$$
\calV\to \calM\to \iota^*\calM^*\to \iota^*\calV^*
$$
is an isomorphism when restricted to $S_q$. Note also that the bundle $\iota^*\calV^*$ is again holomorphic and
$\deg\calV=-\deg\iota^*\calV^*$.

Assume now that $\calM$ is not semi-stable. Then there exists a subbundle $\calV\subset\calM$ of degree $d>0$. But then
$\iota^*\calV^*$ has degree $-d<0$. This implies that for any map $f:\calV\to \iota^*\calV^*$ there exists a non-trivial subbundle
$\calV'\subset \calV$ on which $f$ vanishes. Hence $f$ cannot be injective when restricted to $S_q$. Therefore $\calM$ is semi-stable.

We now want to prove that $\calM$ is split. For this it is enough to prove (in view of \reft{semi-stable}(3))
that if $\calV\subset \calM$ is a subbundle of degree 0, then $\calM$ is isomorphic to $\calV\oplus\calM/\calV$. Let $\calW$ denote the kernel of the natural map $\calM^*\to\calV^*$ (thus $\calW$ is an anti-holomorphic vector bundle).
Since $\calM$ is isomorphic to $\iota^*\calM^*$, the bundle $\iota^*\calW$ is naturally a subbundle of $\calM$ (note that since $\calW$ is anti-holomorphic, the bundle $\iota^*\calW$ has a natural holomorphic structure).
Moreover, it is clear that the corresponding hermitian pairing between $\calV|_{S_q}$ and $(\iota^*\calW)_{S_q}$ is zero. Since this pairing on $\calM|_{S_q}$ is positive-definite, it follows that the fibers of $\calV$ and $\iota^*\calW$ have trivial intersection at any point of $S_q$. This implies
that $\calV\oplus \iota^*\calW$ is a subsheaf of $\calM$
\footnote{Here we identify vector bundles on $E_q$ with their (locally free) sheaves of sections}.
But the degree of $\calM$  (which is equal to 0) is equal to the degree of $\calV\oplus \iota^*\calW$ (which is equal to $\deg\calV+\deg\calW=0+0=0$), hence $\calM$ is isomorphic to $\calV\oplus \iota^*\calW$.

The above analysis shows (still for $G=GL(n)$) that $\eta(\hatg)$ is conjugate to an element of $T$. Let us show that this element is necessarily in $\hatA$ (and thus it can be chosen to lie in $\hatA^+$). We already know that $\calV$ is isomorphic to a direct sum of line bundles
$\calL_1,\cdots,\calL_n$ of degree zero. Hence

1) Every $\calL_i$ is isomorphic to $\iota^*\calL^*$

2) The above isomorphism is positive-definite on $S_q$.

\noindent
Let us assume that $\calL_i=\calL_{a_i}$ where $a_i\in \CC^*$. Then 1) implies that $a_i$ is real and 2) implies that $a_i>0$.
Hence $\eta(\hatg)$ is conjugate to $(a_1,\cdots,a_n,q^2)\in \hatA$.

Let us go back to arbitrary $G$ now. Choosing a faithful unitary representation $(\pi,V)$ of $G$ and using \reft{semi-stable} we
see that the bundle $\calP_{\eta(\hatg)}$ is semi-stable and split. Hence $\eta(\hatg)$ is conjugate to a point in $(a,q^2)\in T\x\CC^*$. Moreover, $(\pi(a),q^2)$ is $*$-invariant, hence $(a,q^2)\in \hatA$.

\ssec{}{Proof of \reft{main}}The above results already show that $\calG_{hol}^+=\calK_{hol}\hatA^+\calK_{hol}$.
Hence, to finish the proof, it is enough to show that $\calK_{hol}\hata_1\calK_{hol}\neq \calK_{hol}\hata_2\calK_{hol}$
if $\hata_1,\hata_2\in\hatA^+$ and $\hata_1\neq \hata_2$. Note that the map $\eta$ is injective on $\hatA^+$. Thus
$\eta(\hata_1)\neq \eta(\hata_2)$. It is well-known that two elements of $\hatA^+$ are  conjugate if and only if they are equal,
hence the statement follows in view of \refl{k-to-g}.


\end{document}